\documentclass{article}

\usepackage{amssymb,amsmath,latexsym,amsthm,eucal}

\usepackage{enumerate}

\newfont{\Bb}{msbm10 scaled\magstep0}
\newfont{\Bbl}{msbm10 scaled\magstep1}
\newfont{\Bbs}{msbm10 scaled 800}

\newcommand{\Z}{\mbox{\Bb{Z}}}
\newcommand{\sZ}{\mbox{\Bbs{Z}}}

\newcommand{\Q}{\mbox{\Bb{Q}}}

\newcommand{\Zp}{\mbox{\Bb{Z}}_p}
\newcommand{\Qp}{\mbox{\Bb{Q}}_p}
\newcommand{\Cp}{\mbox{\Bb{C}}_p}

\newcommand{\W}{B}
\newcommand{\K}{K}

\newcommand{\ord}{{\rm ord}}

\newcommand{\cS}{{\mathbf S}}

\newcommand{\cM}{{\mathbf M}}
\newcommand{\cW}{{\mathbf W}}

\newcommand{\pM}{M}

\newcommand{\Hf}{{\mathbf f}}
\newcommand{\Hg}{{\mathbf g}}
\newcommand{\Hh}{{\mathbf h}}

\newcommand{\RL}{{\mathcal L}}
\newcommand{\E}{{\mathcal E}}
\newcommand{\AJ}{{\rm AJ}}

\newcommand{\f}{\mbox{\Bb F}}

\newenvironment{exafont}{\begin{bf}}{\end{bf}}
\newenvironment{example}{\vspace{0.3cm}\par\noindent\refstepcounter{theorem}\begin{exafont}Example
    \thetheorem\end{exafont}\hspace{\labelsep}}{\vspace{0.3cm}\par}
\newenvironment{algorithm}{\vspace{0.3cm}\par\noindent\refstepcounter{theorem}\begin{exafont}Algorithm
    \thetheorem\end{exafont}\hspace{\labelsep}}{\vspace{0.3cm}\par}
\newenvironment{note}{\vspace{0.3cm}\par\noindent\refstepcounter{theorem}\begin{exafont}Note
    \thetheorem\end{exafont}\hspace{\labelsep}}{\vspace{0.3cm}\par}

\title{Efficient computation of \\ Rankin $p$-adic L-functions}

\author{Alan G.B. Lauder\footnote{Mathematical Institute, 24-29 St Giles, Oxford, United Kingdom. This work is supported in part by a grant from the European
Research Council (204083).} 
\vspace{0.1cm}
\\ Accepted for publication in ``Computation with Modular Forms'' \\ 
Boeckle, G. and Wiese, G. (eds.), Springer Verlag.}

\begin{document}

\maketitle

\begin{abstract}
We present an efficient algorithm for computing certain special
values of Rankin triple product $p$-adic L-functions and give an application of this to the explicit construction of rational points on elliptic curves.
\end{abstract}

\section{Introduction}

The purpose of this paper is to describe an efficient algorithm for computing certain special values of Rankin
triple product $p$-adic L-functions. These special values are $p$-adic numbers and our algorithm
computes them in polynomial-time in the desired precision. This improves on existing algorithms which require exponential
time in the desired precision. Our method has the pleasant feature of also being applicable to Rankin double product $p$-adic L-functions,
and working equally well in weight one as compared to higher weights (Sections \ref{Sec-Weight1} and \ref{Sec-SingleDouble}). We hope it will usefully complement the powerful methods
based upon overconvergent modular symbols for computing $p$-adic L-functions \cite{PS}, which the author understands are less readily adaptable to higher product $p$-adic L-functions.

We describe an application of our algorithm to the efficient construction of rational points on elliptic curves over $\Q$. The curves we consider all
have rank one and relatively small conductor, and so this application does not yield any ``new'' points. However, the constructions give experimental
verification both of the correctness of the implementation of our algorithm, and various sophisticated and new conjectural constructions of rational points on elliptic curves.
Even in the rank one setting these constructions are of interest; for instance, they allow one to carry out by $p$-adic means the complex analytic
calculations in \cite{DDLR} (see Example \ref{Ex-57a}), and in fact
 $p$-adically interpolate points found using much older but not well-understood methods.
A different and enticing application of our algorithm is to the experimental study of conjectural constructions of ``new'' points on elliptic curves over certain number fields using
weight one modular forms \cite{DLR}. In this paper though we shall not address experimentally the 
calculation of these (Stark-Heegner) points attached to weight one forms or the $p$-adic interpolation of points. We
plan to return to these questions in future joint work.

All of the applications of our algorithm are based upon ideas of Darmon and Rotger \cite{DR-Arizona,DR-GrossZagier}. In particular, Darmon encouraged the author 
to try to apply the method for computing with overconvergent modular forms in \cite{AL} to Rankin $p$-adic L-functions, and gave him
invaluable help during the implementation of the algorithm and preparation of this paper. Much of the work behind this paper was in making
the methods in \cite{AL} sufficiently fast in practice to turn a theoretical algorithm (for higher level) into one useful for experimental mathematics.

In writing this paper the author had the choice between trying to give a comprehensive background to the theory necessary to define
Rankin $p$-adic L-functions and present the work of Darmon and Rotger, or distilling just enough to describe his contribution. He chose the
latter, since the long introduction to \cite{DR-GrossZagier} is already very clear but incompressible. This introduction should be read in parallel to
our brief (and simplified) description below by anyone wishing
to get a deeper understanding of the significance of the algorithm in our paper. The reader should also refer to that source for
definitions of any unfamiliar terms below. (All the definitions we shall really need are gathered in Sections \ref{Sec-Background} and
\ref{Sec-RLdef}.)

Let $f,g,h$ be newforms of weights $k,l,m \geq 2$, primitive characters $\chi_f,\chi_g,\chi_h$ with
$\chi_f \chi_g \chi_h = 1$, and level $N$. Assume that the Heegner Hypothesis H \cite[Section 1]{DR-GrossZagier} is satisfied. 
Let $p$ be a prime not dividing the level $N$, and fix an embedding of $\bar{\Q}$ into $\Cp$, the completion of an algebraic closure of the field
of $p$-adic numbers $\Qp$. Assume $f,g$ and $h$ are ordinary at $p$.  Let $\Hf,\Hg$ and $\Hh$ be the (unique) Hida families of (overconvergent)
$p$-adic modular forms passing through $f,g$ and $h$. 
The Rankin $p$-adic L-function $\RL_p^{f}(\Hf,\Hg,\Hh)$ associates to each
triple of weights $(x,y,z)$ in (a suitable subset of) $\Z_{\geq 2}^3$ a $p$-adic number $\RL_p^{f}(\Hf,\Hg,\Hh)(x,y,z) \in \Cp$.
It has a defining interpolation property over
a certain set $\Sigma_{f}$ of unbalanced weights, relating it to the special value of the classical (Garrett-)Rankin triple product L-function
at its central critical point. 
(Weights $(x,y,z)$ are balanced if the largest is strictly
smaller than the sum of the other two, and otherwise unbalanced.)
The theorem of Darmon and Rotger \cite[Theorem 1.3]{DR-GrossZagier} equates its value at the {\it balanced} weights to an explicit algebraic number times 
the $p$-adic Abel-Jacobi map of a certain cycle on a product of Kuga-Sato varieties evaluated at a particular differential form.
At balanced weights $(x,y,z)$ for reasons of sign the classical Rankin triple product L-function vanishes at its central critical point, 
and so the special value $\RL_p^{f}(\Hf,\Hg,\Hh)(x,y,z)$ is thought of as some kind of first derivative. (Darmon and Rotger actually construct in addition
$\RL_p^{g}(\Hf,\Hg,\Hh)(x,y,z)$ and $\RL_p^{h}(\Hf,\Hg,\Hh)(x,y,z)$ but we only consider $\RL_p^{f}(\Hf,\Hg,\Hh)(x,y,z)$ and shall omit from here-on the
superscript $f$.)

In this paper we present an algorithm for computing $\RL_p(\Hf,\Hg,\Hh)(x,y,z) \in \Cp$ for balanced weights $(x,y,z)$
to a given $p$-adic precision in polynomial-time in the precision, provided $p \geq 5$ and {\it under the following assumption on the weights}. Let
us specialise the Hida families back to the original weights $(k,l,m)$ to recover the newforms $f$, $g$ and $h$, and assume that
$(k,l,m)$ is a balanced triple. (This is only a notational simplification --- we are after all really interested in our original newforms, the Hida
families being introduced just to define the interpolation properties of the L-function.) Our algorithm requires that $k = l - m + 2$. This
is enough for all our present and immediately envisaged arithmetic applications.

The problem which makes finding special values of Rankin triple product $p$-adic L-functions challenging is that of
computing ordinary projections of $p$-adic modular forms. That is, in the definition of $\RL_p(\Hf,\Hg,\Hh)(k,l,m)$ one
encounters a $p$-adic modular form ``$d^{-(1+t)}(g^{[p]}) \times h$'' which is not classical, and then has to compute
its ordinary projection to some precision. Since this form is not classical any straightforward approach to this has exponential-time
in the desired $p$-adic precision; for example, by iterating the Atkin operator on $q$-expansions or on some suitable space of
classical modular forms (as the latter necessarily has exponential dimension in the required precision, by consideration of weights cf.
\cite[Proposition I.2.12 ii.]{FG}).

Our solution lies in the fact that ``$d^{-(1+t)}(g^{[p]}) \times h$'' is nearly overconvergent \cite[Section 2.5]{DR-GrossZagier}. 
More precisely, our assumption ($k = l - m + 2$) on the weights
is exactly that which ensures it is overconvergent, and so the methods we developed for computing with such forms in \cite{AL} can be applied.
We expect that our methods can be generalised to handle nearly overconvergent modular forms (using their explicit description in \cite{CGJ})
and thus compute Rankin triple product $p$-adic
L-functions at {\it any} balanced point $(x,y,z)$, but  we have not carried out any detailed work in this direction. The main result of our paper is
really the algorithm (and its refinements) in Section \ref{Sec-ProjAlg} for computing the ordinary projection of certain overconvergent modular forms and
in addition the ordinary subspace. (We give a full and rigorous analysis of this algorithm, but not of two aspects of our overall algorithm for computing Rankin triple product $p$-adic L-functions. These are of minor practical importance, see Note \ref{Note-RankinL} (\ref{Note-NotProved}) and (\ref{Note-HeckeOps}), but difficult to
analyse.)

Regarding arithmetic applications, the most immediate is the following one deduced by Darmon from \cite[Theorem 1.3]{DR-GrossZagier} and \cite[Lemma 2.4]{DRS}, 
in a personal communication.  Assume that $f$ and $g$ are newforms of weight $2$ and trivial character, and that $f$ has rational Fourier coefficients.
Let $E_f$ denote the elliptic curve over $\Q$ associated to $f$, and $\log_{E_f}: E_f(\Qp) \rightarrow \Qp$ be the
formal $p$-adic logarithm map. Then there exists a point $P_g \in E_f (\Q)$ and a computable positive integer $d_g$
such that 
\begin{equation}\label{Eqn-Pg}
\log_{E_f}(P_g) =  2 d_g \frac{\E_0(g) \E_1(g)} {\E(g,f,g)} \RL_p(\Hg,\Hf,\Hg)(2,2,2).
\end{equation}
Here $\E(g,f,g)/\E_0(g) \E_1(g)$ is the explicit non-zero algebraic (in fact quadratic) number which occurs in the Darmon-Rotger formula \cite[Theorem 1.3]{DR-GrossZagier} --- it depends only
upon the $p$th coefficients in the $q$-expansions of $f$ and $g$. Thus if $\RL_p(\Hg,\Hf,\Hg)(2,2,2)$ is non-vanishing one can recover
a point of infinite order on $E_f(\Q)$. 
(The integer $d_g$ is that which appears, in different notation ``$d_T$'' for ``$T := T_g$'', in \cite[Remark 3.1.3]{DDLR}.)
The point $P_g$ is closely related to classically constructed points (``Zhang points'').
We give an example of this application (Example \ref{Ex-57a}), and understand it will be worked out in detail in the forthcoming
Ph.D. thesis of Michael Daub \cite{Daub}. In addition, we also present a number of variations of this application which suggest generalisations of the
different underlying theoretical constructions and also illustrate our algorithm (Section \ref{Sec-Examples}).

The paper is organised in a simple manner, Section \ref{Sec-Algorithm} containing the theoretical background and algorithms, and Section \ref{Sec-Examples} our illustrative computations.\\

{\it Acknowledgements:} This paper would have been neither started nor finished without the constant help and encouragement of Henri Darmon.
It is a pleasure to thank him for this, and to thank also David Loeffler, Victor Rotger and Andrew Wiles for enlightening discussions, and
the anonymous referee for many useful comments.

\section{The Algorithm}\label{Sec-Algorithm}

In this section we present our algorithm for computing the ordinary projection of overconvergent modular forms and certain
special values of Rankin triple product $p$-adic L-functions.

\subsection{Theoretical background}\label{Sec-Background}

We first gather some background material on overconvergent modular forms and the ordinary subspace.

\subsubsection{Katz expansions of overconvergent modular forms}\label{Sec-KatzExp}

Let $N$ be a positive integer, and $p \geq 5$ be a prime not dividing $N$. Let $\chi: (\Z/N\Z)^* \rightarrow \Zp^*$ be a Dirichlet character
with image in $\Zp^*$. The condition that $\chi$ has image in $\Zp^*$ is partly for notational convenience, but see also Note \ref{Note-Minor} (\ref{Note-Fast}).

For each integer $k$ let $\cM_k (N,\chi,\Zp)$ denote the space of classical modular forms for $\Gamma_1(N)$ with character $\chi$
whose $q$-expansions at infinity have coefficients in $\Zp$. This is a free $\Zp$-module of finite rank. 
Let $E_{p-1}$ be the classical Eisenstein series of weight $p-1$ and level $1$ normalised to have constant term $1$.
For each
integer $i > 0$, one may choose a free $\Zp$-module $\cW_i(N,\chi,\Zp)$ of 
$\cM_{k + i(p-1)}(N,\chi,\Zp)$ such that
\[ \cM_{k + i(p-1)}(N,\chi,\Zp) = E_{p-1} \cdot \cM_{k + (i-1)(p-1)}(N,\chi,\Zp) \oplus \cW_i(N,\chi,\Zp).\]
(This choice is not canonical cf. \cite[Page 105]{Katz}.) Define $\cW_0(N,\chi,\Zp) := \cM_k(N,\chi,\Zp)$. 
Let $\K$ be a finite extension of $\Qp$ with ring of integers $\W$.
Define $\cW_i(N,\chi,\W) := \cW_i(N,\chi,\Zp) 
\otimes_{{\sZ}_p} \W$. For $r \in \W$
the space $\pM_k(N,\chi,\W;r)$ of $r$-overconvergent modular forms is by (our) definition the
space of all ``Katz expansions'' of the form
\[ f = \sum_{i = 0}^\infty r^i \frac{b_i}{E_{p-1}^i},\quad b_i \in \cW_i(N,\chi,\W),\quad \lim_{i \rightarrow \infty} b_i = 0\]
where $b_i \rightarrow 0$ as $i \rightarrow \infty$ means the 
$q$-expansions of $b_i$ are more and more divisible by $p$ as $i$ goes to infinity, see \cite[Proposition 2.6.2]{Katz}. 
We define $\pM_k(N,\chi,\K;r) := \pM_k(N,\chi,\W;r) \otimes_{\W} \K$, a $p$-adic Banach space.

The element $r \in B$ plays a purely auxiliary role, determining the inner radius $p^{-\ord_p(r)}$ of the annuli of overconvergence into the supersingular
locus. (Here $\ord_p(\cdot)$ is the $p$-adic valuation normalised with $\ord_p(p) = 1$.) From a computational point of view it is more convenient for each rational number $\alpha > 0$ to consider series of the form
\[ f = \sum_{i = 0}^\infty p^{\lfloor \alpha i \rfloor} \frac{b_i}{E_{p-1}^i},\quad b_i \in \cW_i(N,\chi,\Zp).\]
We just write $\pM_k(N,\chi,\Zp,\alpha)$ for the space of all such elements and call it again the space of $\alpha$-overconvergent modular forms as no
confusion is likely to arise. The space of {\it overconvergent modular forms} $\pM_k(N,\chi,\Zp)$ is the union $\cup_{\alpha > 0} \pM_k(N,\chi,\Zp,\alpha)$.
In everything just defined we may also just forget the character $\chi$ and consider the space $\pM_k(N,\Zp)$ of 
overconvergent modular forms for $\Gamma_1(N)$ itself.
 
\subsubsection{The ordinary subspace}
 
Any overconvergent modular form $f \in \pM_k(N,\Zp)$ is also a $p$-adic modular form \cite[Section 1.4(b)]{JPS} and has a $q$-expansion, and we define the {\it ordinary projection} in the
usual way as $e_{ord}(f)  := \lim_{n \rightarrow \infty} U_p^{n!}(f)$, where $U_p$ is the Atkin operator on $q$-expansions, i.e., $U_p: \sum_n a_n q^n \mapsto
\sum_n a_{np} q^n$. When $k \geq 2$
the image of $e_{ord}$ on $p$-adic modular forms of level $N$ over $\Zp$ is equal to its image on the space of classical modular forms $\cM_k(\Gamma_1(N) 
\cap \Gamma_0(p),\Zp)$ of level $Np$ with trivial character at $p$, see e.g. \cite[Theorem 6.1]{Cole1} or for a precise statement (when $k \geq 3$) \cite[Theorem II.4.3 ii]{FG}. We have for each $\nu \geq 1$ an embedding
\begin{equation}\label{Eqn-nu}
\cM_k(\Gamma_1(N) \cap \Gamma_0(p^\nu),\W) \hookrightarrow \pM_k(N,\W;r)
\end{equation}
for any $r \in \W$ with $\ord_p(r) < 1/p^{\nu - 2}(p+1)$, see (at least for $N \geq 3$) \cite[Corollary II.2.8]{FG}, and also \cite[Page 25]{CGJ}. Thus taking $\nu = 1$ here, one observes for $k \geq 2$
that the image of $e_{ord}$ on $p$-adic modular forms of level $N$ over $\Zp$  is equal (after base change to $\W$) to its image on $\pM_k(N,\W;r)$ for any 
$r \in \W$ with $\ord_p(r) < p/(p+1)$. 
We shall define the $p$-adic {\it ordinary subspace} over $\Zp$ in level $N$, character $\chi$ and weight $k$ to be the image under
$e_{ord}$ of $\pM_k(N,\chi,\Zp,\frac{1}{p+1})$. 
(We make this definition
since this is precisely the space computed by Algorithm \ref{Alg-A}.
For weight $k \geq 2$ this is equivalent to the usual definition as the image of
$p$-adic modular forms under $e_{ord}$, by our preceding observation. The definition should also be equivalent for general weight (certainly
over $\K$) since the ordinary subspace over $\K$ can be described as the space of
overconvergent (generalised) eigenforms of slope zero \cite[Page 59]{FG}, and (generalised) eigenforms of finite slope are
$r$-overconvergent for any $r$ with $\ord_p(r) < p/(p+1)$ \cite[Page 25]{CGJ}.)

\subsection{Projection of overconvergent forms}\label{Sec-ProjAlg}

Underlying our algorithm for computing Rankin $p$-adic L-functions is an algorithm for computing ordinary projections of overconvergent
modular forms and also a basis for the ordinary subspace. It is an extension of  \cite[Algorithm 2.1]{AL}. 

\subsubsection{The basic algorithm}

We first present the basic algorithm, before discussing the steps in more detail and giving some practical refinements. Here
the notation and assumptions are as in Section \ref{Sec-KatzExp}. (We apologise that the
notation ``$m$'' for the $p$-adic precision gives a clash with that used for a weight in the introduction and later, but we
wished to follow closely that in \cite{AL}.)

\begin{algorithm}
\label{Alg-A}
{\it Given an element $H \in \pM_k(N,\chi,\Zp,\frac{1}{p+1})$ where $0 \leq k < p-1$ and integer $m \geq 1$, this algorithm computes
the image in $R:= \Z[[q]]/(p^m,q^{s(m,p)})$ of the ordinary projection $e_{ord}(H)$ and in addition the image in $R$ of an echelonised basis for the ordinary subspace. (Here $s(m,p)$ is some explicit function of $m$ and $p$ defined during the algorithm.)}

\begin{enumerate}[(1)]
\item{[Dimensions]\,
Write $k_0 := k$. Compute $n:= \lfloor \frac{p+1}{p-1}(m+1) \rfloor$. For $i = 0,1,\dots,n$ compute $d_i$, the dimension
of the space of classical modular forms of level $N$ character $\chi$ and weight $k_0 + i(p-1)$. Compute $m_i := 
d_i - d_{i-1}$, for $i \geq 1$,\,$m_0 := d_0$,\, and $\ell := m_0 + m_1 + \cdots + m_n = d_n$. Compute
working precision $m^\prime := m + \lceil \frac{n}{p+1} \rceil$.
Compute $\ell^\prime \geq \ell$, the Sturm bound for the space of classical modular forms of level $N$, character $\chi$ and weight
$k_0 + (p-1)n$.}
\item{[Complementary spaces]\, For each $0 \leq i \leq n$ compute a row-reduced basis $W_i$ of $q$-expansions in
$\Z[[q]]/(p^{m^\prime},q^{\ell^\prime p})$ for some choice of the complementary space $\cW_i(N,\chi,\Zp)$.
}
\item{[Katz expansions]\, Compute the $q$-expansion in $\Z[[q]]/(p^{m^\prime},q^{\ell^\prime p})$ of the Eisenstein series $E_{p-1}(q)$. For  each $0 \leq i \leq n$, let $b_{i,1},\dots,b_{i,m_i}$ denote the elements 
in $W_i$. Compute the ``Katz basis'' elements $e_{i,s} := p^{\lfloor \frac{i}{p+1} \rfloor} E_{p-1}^{-i} b_{i,s}$ in $\Z[[q]]/(p^{m^\prime}, q^{\ell^\prime p})$.}
\item{[Atkin operator]\, For each $0 \leq i \leq n$ and $1 \leq s \leq m_i$ compute $t_{i,s} := U_p(u_{i,s})$ in $\Z[[q]]/(p^{m^\prime},q^{\ell^\prime})$, where $U_p$ is the Atkin operator on $q$-expansions and $u_{i,s} := e_{i,s}$.}
\item{[Atkin matrix]\, Compute
$T$, the $\ell \times \ell^\prime$ matrix over $\Z/(p^{m^\prime})$ whose entries are the coefficients
in the $q$-expansions modulo $q^{\ell^\prime}$ of the $\ell$ elements $t_{i,s}$. Compute $E$, the
$\ell \times \ell^\prime$ 
matrix over $\Z/(p^{m^\prime})$ whose entries are the coefficients in the $q$-expansions
modulo $q^{\ell^\prime}$ of the $\ell$ elements $e_{i,s}$. Use linear algebra over $\Z/(p^{m^\prime})$ to compute the matrix
$A^\prime$ over $\Z/(p^{m^\prime})$ such that $T = A^\prime E$. Let $A$ be the ``Atkin matrix'' over $\Z/(p^m)$ obtained by reducing entries
in $A^\prime$ modulo $p^m$.}
\item{[Two-stage projection]
Compute the image $H \in \Z[[q]]/(p^{m^\prime},q^{\ell^\prime p})$.
\begin{enumerate}
\item{[Improve overconvergence]
Compute $U_p(H) \in \Z[[q]]/(p^{m^\prime},q^{\ell^\prime})$ and find coefficients $\alpha_{i,s} \in \Z/(p^m)$ such that $U_p(H) \equiv \sum_{i,s} \alpha_{i,s} e_{i,s} \bmod{(p^{m},q^{\ell^\prime})}$.}
\item{[Projection via Katz expansion]
Compute
a positive integer $f$ such that all the unit roots of the reverse characteristic polynomial of $A$ lie in some extension of $\Zp$ with
residue class field of degree $f$ over $\f_p$. Compute $A^{r-1}$  for $r:= (p^f - 1)p^m$ using fast exponentiation. Compute
$\gamma := \alpha A^{r-1}$ where $\alpha$ is the row vector $(\alpha_{i,s})$. 
Write $\gamma=  (\gamma_{i,s})$ and return the ordinary projection $e_{ord}(H) = \sum_{i,s} \gamma_{i,s} e_{i,s} \in \Z[[q]]/(p^m,q^{s(m,p)})$ where $s(m,p) := \ell^\prime p$.}
\end{enumerate}
}

\item{[Ordinary subspace]
Compute $A^r = A^{r-1} A$ and let $\{(B_{i,s})\}$ be the set of non-zero rows in the echelon form $B$ of the matrix $A^r$. Return 
$\sum_{i,s} B_{i,s} e_{i,s}  \in \Z[[q]]/(p^m,q^{s(m,p)})$ for each non-zero row $(B_{i,s})$,  the image of a basis for the ordinary subspace.}
\end{enumerate}
\end{algorithm}

In this algorithm we assume that the $q$-expansion of the input modular form $H$ can be computed in polynomial-time in $N,p$ and any desired $p$-adic and $q$-adic precisions.
Regarding the complexity of the whole algorithm, we just refer the reader to the analysis of Steps 1-5 in \cite[Sections 3.2.2, 3.3.1]{AL}, and
observe that Steps 6 and 7 can be carried out using standard methods in linear algebra. In particular, the algorithm
is certainly polynomial time in $N,p$ and $m$.

\subsubsection{Proof of correctness}

The analysis of the correctness of the algorithm is very similar to that in \cite[Section 3.2.1]{AL}. The essential idea is the following. One considers an infinite square matrix
for the Atkin $U_p$ operator on the space of $\frac{1}{p+1}$-overconvergent modular forms w.r.t. some choice of Katz basis. Reducing this
(assumed integral, see Note \ref{Note-Minor} (\ref{Note-A}))
matrix modulo $p^m$, it vanishes
except for an $\infty \times \ell$ strip down the lefthand side. The matrix $A$ modulo $p^m$ we compute is the $\ell \times  \ell$ matrix which occurs in the top lefthand corner, for our
choice of basis (this is proved in \cite[Section 3]{AL}). 
We would like to iterate the infinite matrix on the infinite row vector representing an overconvergent modular form $H$. When
$H \in \pM_k(N,\chi,\Zp,\frac{p}{p+1})$ we notice that the coefficients in the infinite vector representing $H$ w.r.t. our Katz basis decay
$p$-adically (since $p/(p+1) > 1/(p+1)$) and in fact vanish modulo $p^m$, except for
the first $\ell$ elements (see the final paragraph in \cite[Section 3.4.2]{AL}).
Hence we can iterate $U_p$ on $H$ by iterating the finite matrix $A$ on a finite vector of length $\ell$. (The actual power $r$ is chosen to
ensure that we iterate sufficiently often to obtain the correct answer modulo $p^m$.)
In our application to Rankin $p$-adic L-functions we will find that in fact $H \in \pM_k(N,\chi,\Zp,\frac{1}{p+1})$. 
Hence the preliminary Step 6 (a) is to apply the $U_p$ operator once to
$q$-expansions to improve overconvergence by a factor $p$, see \cite[Equation (2)]{AL} and Note \ref{Note-Minor} (\ref{Note-A}). (There is a loss of precision of $m^\prime - m$ when one writes
$U_p(H)$ as a Katz expansion, cf. the last paragraph of \cite[Section 3.2.1]{AL} where a similar loss occurs during the computation of the matrix $A$.)
Observe that this preliminary step is harmless, since we need to compute the elements in our Katz basis to the higher precision
modulo $q^{\ell^\prime p}$ anyway. (To make the above argument completely rigorous one fusses over the minor difference between
$r$-overconvergent for all $\ord_p(r) < p/(p+1)$, and $\frac{p}{p+1}$-overconvergent, as in \cite[Section 3.2.1]{AL}.)

\begin{note}\label{Note-Minor}
We make some minor comments on the algorithm.
\begin{enumerate}[(1)]

\item{
For weight $k \geq 2$ the ordinary subspace can be computed instead using classical methods; however, our algorithm
is the only ``polynomial-time'' method known to the author for computing this subspace in weight $k \leq 1$. 
}\label{Note-k1}

\item{
We assume that the smallest non-zero slope $s_0$ of (the Newton polygon of) the reverse characteristic polynomial of $A$ is 
such that $\lceil m/s_0 \rceil \leq (p^f - 1)p^m$. This is 
reasonable as the smallest non-zero slope which has ever been experimentally observed is $1/2$. (One could of course compute $s_0$ and adjust
$r$ accordingly to remove this assumption.) The integer $f$ can be 
easily computed by reducing the matrix $A$ modulo $p$. The exponent $m$ rather than $m-1$ in the definition of $r$ accounts for the possibility
that the unit roots may lie in ramified extensions. (So $u^r \equiv 1 \bmod{p^m}$ for each unit root $u$, and $u^r \equiv 0 \bmod{p^m}$ for all other roots $u$ of the reverse characteristic polynomial.)
}\label{Note-s0}

\item{
The correctness of the algorithm relies on the assumption that we can solve $T = A^\prime E$ for a $p$-adically integral
matrix $A^\prime$, although the theory only guarantees that $p A^\prime$ has integral coefficients, see \cite[Note 3.2, Section 3.2.1]{AL} and 
also \cite[II.3]{FG}, \cite[Section 3.11]{Katz}.  One could modify the algorithm (or rather the refined version in Section \ref{Sec-S6}) to remove this
assumption; however, in practice the author has never encountered a situation in which the matrix
$A^\prime$ fails to have $p$-adically integral coefficients.
}\label{Note-A}

\item{
The assumption that $\Z[\chi]$ embeds in $\Zp$ allows one to exploit fast algorithms for matrix and polynomial arithmetic over rings of the form
$\Zp/(p^m) \cong \Z/(p^m)$ which are integrated into the systems {\sc Magma} and {\sc Sage}. 
The algorithm works perfectly well in principle without this assumption, but it will be much more difficult to get a comparably fast implementation.
}\label{Note-Fast}

\item{
The hypothesis $0 \leq k < p-1$ can be removed as follows. In Step 1 write $k:= k_0 + j(p-1)$ where $0 \leq k_0 < p-1$. In Step 4 compute
$G := E_{p-1}(q)/E_{p-1}(q^p)$ and $G^j \in \Z[[q]]/(p^{m^\prime},q^{\ell^\prime p})$, and let $u_{i,s} := G^j e_{i,s}$. The matrix $A$ computed in Step 5 is
then for the ``twisted'' Atkin operator $U_p \circ G^j$. After Step 6(a) multiply the $q$-expansion of $U_p(H)$ by $E_{p-1}^{-j}$ and in Step 6(b)
multiply the $q$-expansion $\sum_{i,s} \gamma_{i,s} e_{i,s}$ by $E_{p-1}^j$ and return this product as $e_{ord}(H)$.  
In Step 7 multiply each $q$-expansion $\sum_{i,s} b_{i,s} e_{i,s}$ by $E_{p-1}^j$ to get the basis for ordinary space. For $j \geq 1$, to ensure
$E_{p-1}^{-j} U_{p}(H)$ lies in the correct space one should multiply it by $p^{\lceil \frac{j}{p+1} \rceil}$, and so the final
answer will only be correct modulo $p^{m - \lceil \frac{j}{p+1} \rceil}$. (One could of course also just run the algorithm without twisting, but then
the auxiliary parameters $n,\ell, m^\prime$ etc would have to be worked out afresh, since the algorithm would no longer be an extension
of \cite[Algorithm 2.1]{AL}.)
}\label{Note-kbig}

\item{
In practice the output $q$-adic precision $s(m,p) = \ell^\prime p $ is always large enough for our needed application to Rankin $p$-adic L-functions. One
can insist though on any output precision $s^\prime \geq \ell^\prime p$ simply by computing the Katz basis elements in Step 3 to that $q$-adic precision.
}\label{Note-qadic}

\end{enumerate}
\end{note}

\subsubsection{Finding complementary spaces in Step 2}\label{Sec-S23}

A key step in the algorithm is the efficient construction in practice of the image in
$ \Z[[q]]/(p^{m^\prime},q^{\ell^\prime p})$ of a basis for some choice of complementary spaces $\cW_i(N,\chi,\Zp)$, for each
$0 \leq i \leq n$. The author's implementation (which is
at present restricted to trivial or quadratic characters $\chi$) is based upon suggestions of David Loeffler and John Voight. The idea is to use the multiplicative structure
on the ring of modular forms. 

One fixes a choice of weight bound ${\mathcal B} \geq 1$ and computes
the image in $\Z[[q]]/(p^{m^\prime},q^{\ell^\prime p})$ of a $\Z[\chi^\prime]$-basis for each of the spaces of classical modular forms $\cM_b  (N,\chi^\prime,\Q(\chi^\prime))$ where $1 \leq b \leq {\mathcal B}$ and $\chi^\prime$ vary over a set of characters which generate a group containing $\chi$. One then reduces these basis elements modulo
$(p,q^{\ell^\prime})$ and for each $0 \leq i \leq n$ looks for random products of these $q$-expansions which generate an $\f_p$-vector space of dimension $d_i$ and have weight
$k_0 + i(p-1)$ and character $\chi$. This is done in a recursive manner. Once one has computed
the required forms in weight $k_0 + (i-1)(p-1)$ one maps then (via the identity map) into weight $k_0 + i (p-1)$ (recall $E_{p-1} \equiv 1 \mod{p}$) and 
generates a further $m_i  = d_i - d_{i-1}$ linearly independent forms in weight $k_0 + i (p-1)$ and character $\chi$. The correct choices of
products, which give forms not in the space already generated,  are ``encoded'' in an appropriate manner; that is, the basis elements for each weight $b$ and character $\chi^\prime$ are stored as an ordered list, and products of them (modulo $(p,q^{\ell^\prime})$) can then be represented
by ``codes'' which give the positions chosen in each list.

Having found these correct choices modulo $(p,q^{\ell^\prime})$ one then repeats the process modulo $(p^{m^\prime},q^{\ell^\prime p})$ to find the complementary spaces $\cW_i(N,\chi,\Zp)$,
but crucially this time using the ``codes'' to only take products of modular forms which give something not in the $\Zp$-span of the forms already computed.
In this way when working to the full precision
one does not waste time computing products of modular forms that lie in the space one has already generated. (It surprised the author to discover that in practice
in some examples, e.g. Example \ref{Ex-469}, one can generate many such ``dud'' forms --- he has no intuition as to why this is the case.)

A good bound to take is ${\mathcal B} := 6$, but one can vary this, playing the time it takes to generate the spaces in low weight off against the time spent looking for
suitable products. This choice of bound fits with some theoretical predictions communicated to the author by David Loeffler.

\subsubsection{A three-stage projection in Step 6}\label{Sec-S6}

In Algorithm \ref{Alg-A} we find $U_p^r(H)$, where $r$ is chosen so that the answer
is correct modulo $p^m$, in two separate stages. First, one computes $U_p(H) \in \pM_k(N,\chi,\Zp,\frac{p}{p+1})$ using
$q$-expansions. Second, one computes $U_p^{r-1}(U_p(H))$ using Katz expansions. However, the matrix $A$ has size
growing linearly with $m$
and so the computation of the high power $A^{r-1}$ becomes a bottleneck as the precision $m$ increases.

A better aproach is to factor the projection map into three parts, as follows. Write $s_0$ for the smallest non-zero slope in the
characteristic series of $A$ (one can safely just set $s_0 := 1/2$).
Computing $A^{\lceil m/s_0 \rceil }$ and writing
its non-zero rows (which are w.r.t. the Katz basis) as $q$-expansions in $\Z[[q]]/(p^m,q^{\ell^\prime p})$ gives
(the image of) a basis for the ordinary subspace. One can now compute a matrix $A_{ord}$ over $\Z/(p^m)$ for the $U_p$ operator
on this basis by explicitly computing with $q$-expansions. This matrix is significantly smaller than $A$ itself, since its dimension
has no dependence on $m$. To project
$H$, one computes as before $U_p(H)$ using $q$-expansions, then $U_p^{\lceil m/s_0 \rceil}(U_p(H))$ via Katz expansions as the product
$\beta := \alpha A^{\lceil m/s_0 \rceil}$. Next, one writes the ``Katz vector'' $\beta$ as the image of a $q$-expansion in
$\Z[[q]]/(p^m,q^{\ell^\prime p})$ and thus as a new vector $\beta^{\prime}$ over $\Z/(p^m)$ in terms of the basis for the ordinary subspace. Finally, 
one computes $U_p^{r - \lceil m/s_0 \rceil - 1}$ on $U_p^{\lceil m/s_0 \rceil}(U_p(H))$ as $\gamma^\prime := \beta^\prime A_{ord}^{r - \lceil m/s_0 \rceil - 1}$
and returns the $q$-expansion associated to $\gamma^\prime$ as the ordinary projection of $H$ modulo $(p^m,q^{s(m,p)})$.

This three-stage projection method also works for $k < 0$ or $k \geq p-1$, but one must take care to twist and un-twist by powers of $E_{p-1}$
at the appropriate times.

\subsubsection{Avoiding weight one forms}\label{Sec-Weight1}

In the case that $k = 1$, one can compute the ordinary projection $e_{ord}(H)$ of the weight one form $H$ without doing {\it any} computations 
in weight one, except computing the $q$-expansion of $H$ itself modulo $(p^{m^\prime}, q^{\ell^\prime p})$. The idea is to use the Eisenstein series
to ``twist'' up to weight $p = 1 + (p-1)$. That is, one proceeds as in Note \ref{Note-Minor} (\ref{Note-kbig}), only writing  $k =  1 = k_0 + j(p-1)$ where now $k_0 := p$ and
$j:=-1$.  In addition, when generating complementary spaces (see Section \ref{Sec-S23}) one only computes bases of classical modular forms in low weights
$2 \leq b \leq {\mathcal B}$.

The author has implemented this variation in both {\sc Magma} and {\sc Sage}, and used it to compute the characteristic series of the Atkin operator on $p$-adic overconvergent
modular forms in weight one (for a quadratic character, and various levels $N$ and primes $p$) without computing the $q$-expansions of any modular forms in weight one.
 
\subsection{Application to $p$-adic L-functions}\label{Sec-RL}

We now describe the application of Algorithm \ref{Alg-A} to the computation of $p$-adic L-functions.

\subsubsection{Definition of Rankin triple product $p$-adic L-functions}\label{Sec-RLdef}

Let $f,g,h$ be newforms of balanced weights $k,l,m \geq 2$, primitive characters $\chi_f,\chi_g,\chi_h$, with $\chi_f \chi_g \chi_h = 1$ and level $N$.
Assume that the Heegner hypothesis H from \cite[Section 1]{DR-GrossZagier}
is satisfied, e.g. $N$ is squarefree and for each prime $\ell$ dividing $N$ the product of the $\ell$th Fourier
coefficients of $f,g$ and $h$ is $-\ell^{(k + l + m - 6)/2}$. 
Write $k = l + m - 2 - 2t$ with $t \geq 0$, which is possible since the sum of the weights must be even. 
We fix an embedding $\bar{\Q} \hookrightarrow \Cp$ and assume $f,g$ and $h$ are ordinary at $p$. That is, the
$p$th coefficient in the $q$-expansion of each is a $p$-adic unit.

Define the map $d = q \frac{d}{dq}$ on $q$-expansions as $d: \sum_{n \geq 0} a_n q^n \mapsto \sum_{n \geq 0} n a_n q^n$.
Then for $s \geq 0$, the map $d^s$ acts on $p$-adic modular forms increasing weights by $2s$ \cite[Th\'eor\`eme 5(a)]{JPS}. For a $p$-adic modular form
$a(q) := \sum_{n \geq 0} a_n q^n$ let $a^{[p]} := \sum_{n \geq 1,\, p \not \;| n} a_n q^n$ denote its $p$-depletion. Then for
$s \geq 1$ the map 
\[ a  \mapsto d^{-s}(a^{[p]}) = \sum_{\stackrel{n = 1}{p \not \;| \,n}}^{\infty} \frac{a_n}{n^s} q^n\]
acts on $p$-adic modular forms shifting weights by $-2s$ \cite[Th\'eor\`eme 5(b)]{JPS}.
So $d^{-(1 + t)}(g^{[p]}) \times h$ is a $p$-adic modular form of weight $l - 2(1 + t) + m = k$ and character $\chi_g \chi_h = \chi_f^{-1}$.

 Let $f^*$ be the dual form to $f$ and
$f^{*(p)}$ be the ordinary $p$-stabilisation of $f^*$, see \cite[Sections 1 and 4]{DR-GrossZagier}. 
So $f^{*(p)}$ is an ordinary eigenform of character $\chi_f^{-1}$.
We define
\[ \RL_p(\Hf,\Hg,\Hh)(k,l,m) := c(f^{*(p)}, e_{ord}(d^{-(1 + t)}(g^{p]}) \times h)) \in \Cp.\]
Here we are assuming the action of the Hecke algebra on the ordinary subspace in weight $k$ is
semisimple (which the author understands is well-known for $N$ squarefree since $k \geq 2$)
and $c(f^{*(p)},\bullet)$ denotes the coefficient of $f^{*(p)}$ when one writes
an ordinary form $\bullet$ as a linear combination of ordinary eigenforms, see \cite[Page 222]{Hida-Book}. (Darmon and
Rotger take a different but equivalent approach, using the Poincar\'e pairing in algebraic de Rham cohomology to extract the coefficient
$\RL_p(\Hf,\Hg,\Hh)(k,l,m)$ \cite[Proposition 4.10]{DR-GrossZagier}.)

\subsubsection{Computation of Rankin triple product $p$-adic L-functions}

We shall now (with another apology) introduce the clashing notation $m$ to refer to the $p$-adic precision, as in Section \ref{Sec-ProjAlg}.
We wish to apply Algorithm \ref{Alg-A} to compute $e_{ord}(H)$ for $H := d^{-(1 + t)}(g^{[p]}) \times h$ modulo $p^m$
(and $q$-adic precision $s(m,p)$) and
also a basis for the ordinary subspace in level $N$, weight $k$ and character $\chi:= \chi_f^{-1}$.
(So we should assume that the image of $\chi$ lies in $\Zp$ and $g$ and $h$ are defined over $\Zp$, but see also Notes \ref{Note-Minor} (\ref{Note-Fast}) and
\ref{Note-RankinL} (\ref{Note-phi}).)
Given these, one can use Hecke operators on the ordinary subspace to extract the coefficient
$c(f^{*(p)},e_{ord}(H))$, see Note \ref{Note-RankinL} (\ref{Note-HeckeOps}), as $f^{*(p)}$ (and $H$) are easy to compute (at least within {\sc Magma} and {\sc Sage} using the algorithms
from \cite{Stein}).

For our projection algorithm to work we require that $H$ is overconvergent (rather than just nearly overconvergent \cite[Section 2.5]{DR-GrossZagier}) and in particular that 
$H \in \pM_k(N,\chi,\Zp,\frac{1}{p+1})$ for $\chi := \chi_f^{-1}$. Overconvergence is guaranteed provided $t = l - 2$, since
$a \mapsto d^{-s}(a^{[p]})$ maps overconvergent forms in weight $1 + s$ to overconvergent forms in weight $1 - s$ \cite[Proposition 4.3]{Cole1}. 
That is, provided
\begin{equation}\label{Eqn-klm}
k = m - l + 2
\end{equation}
we will have that $d^{-(1+t)}(g^{[p]})$, and hence also $H$, is overconvergent.
When this condition is not satisfied our algorithm with fail.

Regarding the precise radius of convergence of $d^{-(1+t)}(g^{[p]})$, Darmon has explained to the author that when our condition (\ref{Eqn-klm}) is met the methods
used by Coleman (the geometric interpretation of the $d$ operator in terms of the Gauss-Manin connection \cite[Section 2.5]{DR-GrossZagier})
show the form $d^{-(1 + t)}(g^{[p]})$ lies in the space $\pM_k(N,\chi_g,\K;r)$ for any 
$r \in \W \subset \Cp$ with $\ord_p(r) < 1/(p+1)$. Let us outline the argument to get an idea why this is true. First, the $p$-depletion 
 $g^{[p]} := (1 - V_p U_p)g$ 
 is a classical modular form for $\Gamma_1(N) \cap \Gamma_0(p^2)$ with trivial character at $p$ and infinite slope. 
 (Here $V_p$ is the one-sided inverse of $U_p$ \cite[Equation (12)]{DR-GrossZagier} and increases the level of $U_p (g)$ by $p$.) 
Hence by (\ref{Eqn-nu}) with $\nu := 2$, $g^{[p]}$ lies in $\pM_{\ell}(N,\chi_g,\W;r)$ for any $r \in \W$ with $\ord_p(r) < 1/(p+1)$.
Next, \cite[Theorem 5.4]{Cole1} gives an explicit relation between the action of powers of the $d$ operator on
spaces of overconvergent modular forms and that of the Gauss-Manin connection of certain de Rham cohomology spaces associated
to rigid analytic modular curves. This relationship associates to $g^{[p]}$ a trivial class in the de Rham cohomology space (the
righthand side of \cite[Equation (34)]{DR-GrossZagier} for ``$r$'' equals $t$ and any ``$\epsilon$'' less than $1/(p+1)$), and hence
one in the image of the Gauss-Manin connection. (The class is trivial because the form has infinite slope, cf. \cite[Lemma 6.3]{Cole1}.)
The Gauss-Manin connection preserves the radius of convergence, and taking the preimage  and untangling
the relationship one finds that $d^{-(1+t)}(g^{[p]})$ is an overconvergent modular form of the same radius of convergence as $g^{[p]}$, i.e.
$d^{-(1+t)}(g^{[p]}) \in \pM_{\ell}(N,\chi_g,\K;r)$ for any $r \in \W$ with $\ord_p(r) < 1/(p+1)$.
Thus multiplying by $h$ (and using (\ref{Eqn-nu}) with $\nu := 1$ to determine the overconvergence of $h$ itself) we find also
\begin{equation}\label{Eqn-dgh}
d^{-(1 + t)}(g^{[p]}) \times h \in \pM_k(N,\chi,\K;r)
\end{equation}
for any $r \in \W$ with $\ord_p(r) < 1/(p+1)$.

\begin{note}\label{Note-RankinL}
\begin{enumerate}[(1)]
\item{
The above argument does not quite show that $H := d^{-(1 + t)}(g^{[p]}) \times h$ lies in $\pM_k(N,\chi,\Zp,\frac{1}{p+1})$
as for this one would need to replace ``$\K$'' by ``$\W$'' in (\ref{Eqn-dgh}).
However, the author just {\it assumed} this was true, and this was not contradicted by our experiments; in particular,
when one could relate the value of the Rankin $p$-adic L-function to the $p$-adic logarithm of a point on an elliptic curve, the relationship held to exactly the 
precision predicted by the algorithm. To be completely rigorous though one would have to carry out a detailed analysis of Darmon's argument and
the constructions used by Coleman (and
one may have to account for some extra logarithmic growth and loss of precision).
}\label{Note-NotProved}
\item{
It is helpful to notice that the map $\phi:(g,h) \mapsto e_{ord}(d^{-(1+t)}(g^{[p]}) \times h)$ is bi-linear in $g$ and $h$. Thus one can
compute $\phi((g,h))$ by first computing it on a product of bases for the spaces $\cS_l(N,\chi_g)$ and $\cS_m(N,\chi_h)$. This is
useful when these spaces are defined over $\Zp$ but the newforms themselves are defined over algebraic number fields which
do not embed in $\Zp$.
}\label{Note-phi}
\item{
The author implemented a number of different approaches to computing $c(f^{*(p)}, e_{ord}(H))$. The most straightforward is
to compute matrices for the Hecke operators $U_\ell$ (for $\ell | Np)$ and $T_\ell$ (otherwise) on the ordinary subspace for many small
$\ell$ by explicitly computing on the $q$-expansion basis for the ordinary subspace output by Algorithm \ref{Alg-A}. One can then try to project onto the ``$f^{*(p)}$-eigenspace''
using any one of these matrices. One difficulty which arises is that congruences between eigenforms may force a small loss of $p$-adic precision during this
step. (Congruences with Eisenstein series can be avoided for $k \geq 2$ by using classical methods to compute a basis for the ordinary cuspidal
subspace, and working with that space instead.) We did not carry out a rigorous analysis of what loss of precision could occur
due to these congruences, but in our examples it was never more than a few $p$-adic coefficients and one could always determine exactly what loss
of precision had occurred after the experiment. The author understands from discussions with Wiles that one should be able to compute an ``upper bound'' on 
the $p$-adic congruences which can occur, and thus on the loss of precision. This bound of course is entirely independent of the precision $m$.
However,  such a calculation is beyond the scope of this paper.
}\label{Note-HeckeOps}
\end{enumerate}
\end{note}

\subsubsection{Single and double product L-functions}\label{Sec-SingleDouble}

The author understands that usual $p$-adic L-functions can be computed using our methods, by substituting Eisenstein series for
newforms in the appropriate places in the triple product L-function, cf. \cite[Section 3]{BD}.  However, he has not looked at this application at all, as the methods based upon overconvergent modular symbols are already very good (for $k \geq 2$) \cite{PS}. One can similarly compute
double product Rankin $p$-adic L-functions using our approach. In particular, we have used our algorithm to compute a (suitably defined) Rankin double product 
$p$-adic L-function special value ``$\RL_p(\Hf,\Hg)(2,1)$'' for $f$ of weight $2$ and $g$ of weight $1$, see the forthcoming \cite{DLR} and also
\cite[Conjecture 10.1]{DR-Arizona}. 

\section{Examples}\label{Sec-Examples}

In this section we shall freely use the notation from \cite[Section 1]{DR-GrossZagier}.
We implemented our basic algorithms in both {\sc Magma} and {\sc Sage}, but focussed our refinements on the former and all the examples
we present here were computed using this package. The running time and space for the examples varied from 
around 100 seconds with 201 MB RAM (Example \ref{Ex-77a}) to around 19000 seconds with 9.7 GB RAM (Example \ref{Ex-43a}) on a
2.93 GHz machine.



All of the examples here are for weights $k,l,m$ with $k = m - l + 2$ and $t = l - 2$, where $t = 0$, i.e., $l = 2$ and $k = m$ (and in fact $f = h$). We implemented our algorithm for arbitrary $t \geq 0$ and
computed $\RL_p(\Hf,\Hg,\Hh)(k,l,m)$ in cases when $t > 0$; however, the author does not know of any geometric constructions of points when $t > 0$ (or even when $f \ne h$) 
and so we do not present these computations here.

We begin with an  example of the explicit construction of rational points mentioned in our introduction, see
Equation (\ref{Eqn-Pg}).

\begin{example}\label{Ex-57a}
Let $E_f\,:\, y^2 + y = x^3 - x^2 - 2x + 2$
be the rank $1$ curve of conductor $57$ with Cremona label ``57a'' associated to the cusp form
\[ f:= q - 2q^2 - q^3 + 2q^4 - 3q^5 + 2q^6 - 5q^7 + q^9 + \cdots .\]
We choose two other newforms of level $57$ (associated to curves of rank zero):
\[
\begin{array}{rcl}
g_1  & := & q + q^2 + q^3 - q^4 - 2q^5 + q^6 - 3q^8 + q^9  - \cdots \\
g_2  & := & q - 2q^2 + q^3 + 2q^4 + q^5 - 2q^6 + 3q^7 + q^9  - \cdots.
\end{array}
\]
Taking $p := 5$ and writing $\Hf,\Hg_1,\Hg_2$ for the Hida families we compute the special values
\[
\begin{array}{rcl}
\RL_5(\Hg_1,\Hf,\Hg_1)(2,2,2) & \equiv &  -260429402433721822483 \bmod{5^{30}}\\
5 \RL_5(\Hg_2,\Hf,\Hg_2)(2,2,2) & \equiv &  -279706401244025789341 \bmod{5^{31}}.
\end{array}
\]
One computes that for each newform $g_i$, if one multiplies the operator of projection onto the
$g_i$-eigenspace by $3$ then one obtains an element in the integral (rather than rational) Hecke algebra.
Thus equation (\ref{Eqn-Pg}) predicts that there exist global points $P_1,P_2 \in E_f(\Q)$ such that
\[  \log_{E_f}(P_i) = 6 \times \frac{\E_0(g_i) \E_1(g_i)}{\E(g_i,f,g_i)} \times \RL_5(\Hg_i,\Hf,\Hg_i) .\]
One finds
\[
\begin{array}{rcl}
\log_{E_f}(P_1) & \equiv  &  37060573996879427247 \times 5 \bmod{5^{30}}\\
\log_{E_f}(P_2) & \equiv & -18578369245374641968 \times 5 \bmod{5^{30}}.
\end{array}
\]
Adapting the method in \cite[Section 2.7]{KP} we recover the points
\[
\begin{array}{rcccl}
P_1 & = & \left(-\frac{1976}{7569},\frac{750007}{658503}\right) & = & -16 P
\end{array}
\]
and $P_2  = (0,1) =  4P$,
where $P := (2,-2)$ is a generator for $E_f(\Q)$.
\end{example}

Next we look at an example where the Darmon-Rotger formula may be applied, but the application to constructing
points has not been fully worked out. (At least, at the time of author's computations --- we understand from a personal communication from
Darmon and Rotger that this has now been done.)

\begin{example}
Let $E_f\,:\, y^2 + xy + y = x^3 - x^2$
be the rank $1$ curve of conductor $53$ with Cremona label ``53a'' associated to the cusp form
\[ f:= q - q^2 - 3 q^3 - q^4 + 3 q^6 - 4 q^7 + 3 q^8 + 6 q^9  + \cdots .\]
There is one newform $g$ of level $53$ and weight $4$ and trivial character with
rational Fourier coefficients:
\[ g := q + q^3 - 8 q^4 - 18 q^5 + 2 q^7 - 26 q^9 + 54 q^{11} + \cdots.\]
Taking $p:=7$ and writing $\Hf$ and $\Hg$ for the Hida families we compute the special value
 \[ \RL_7(\Hg,\Hf,\Hg)(4,2,4) \equiv  -12581507765759084963366603 \bmod{7^{30}}.\]
The Darmon-Rotger formula \cite[Theorem 1.3]{DR-GrossZagier} then predicts that
\[ \AJ_7(\Delta)(\eta_g^{{\rm u-r}} \otimes \omega_f \otimes \omega_g) = \frac{\E_0(g) \E_1(g)}{\E(g,f,g)}
\RL_7(\Hg,\Hf,\Hg)(4,2,4)\] 
and we find that
\[ \AJ_7(\Delta)(\eta_g^{{\rm u-r}} \otimes \omega_f \otimes \omega_g) \equiv 1025211670724558054729221 \times 7  \bmod{7^{30}}.\]
Equation (\ref{Eqn-Pg}) does not apply in this setting, but one can hope that this equals $\log_{E_f}(P)$ for
some point $P  \in E_f(\Q) \otimes \Q$. Exponentiating one finds a point
$\hat{P}  = (x(\hat{P}),y(\hat{P})) \in E_1(\Q_7)$ 
with coordinates $7^2 x(\hat{P}),\,7^3 y(\hat{P})$ modulo $7^{30}$
(where $E := E_f$). We have $|E(\f_7)| = 12$ and translating
$\hat{P}$ by elements $Q \in E(\Q_7)[12]$ we find exactly one rational point, $P = (0,-1)$ (see the method 
in \cite[Section 2.7]{KP}). Thus we have computed
a generator in a rather elaborate manner.
\end{example}

The author also considered again
the curve $E_f$ with Cremona label ``57a'' but took $g$ to be the unique newform of level
$57$ and weight $4$ with trivial character and rational Fourier coefficients, and found that
$\AJ_5(\Delta)(\eta_g^{{\rm u-r}} \otimes \omega_f \otimes \omega_g) \equiv - \frac{15}{13} \log_{E_f}(P) \bmod{5^{31}}$ for
$P := (2,-2)$ a generator of $E_f(\Q)$. (So here $p = k - 1$, and we used the ``twisted'' version of the algorithm
described in Note \ref{Note-Minor} (\ref{Note-kbig}).) 

The next example has a similar flavour but involves cusp forms of odd weight.

\begin{example}\label{Ex-43a}
Let $E_f\,:\, y^2 + y = x^3 + x^2$
be the rank $1$ curve of conductor $43$ with Cremona label ``43a'' associated to the cusp form
\[ f:=  q - 2 q^2 - 2 q^3 + 2 q^4 - 4 q^5 + 4 q^6 + q^9 + 8 q^{10} + 3 q^{11} + \cdots .\]
Let $\chi$ be the Legendre character modulo $43$. Then we find unique newforms
$g \in S_3(43,\chi)$ and $h \in S_5(43,\chi)$ with rational Fourier coefficients:
\[
\begin{array}{rcl}
g & := & q + 4 q^4 + 9 q^9 - 21 q^{11} + \cdots \\
h & := & q + 16 q^4 + 81 q^9 + 199 q^{11} + \cdots.
\end{array}
\]
Taking $p:=11$ and writing $\Hf,\Hg$ and $\Hh$ for the Hida families we compute the special values
\[ 
\begin{array}{rcl}
 \RL_{11}(\Hg,\Hf,\Hg)(3,2,3) & \equiv &-7831319270947510009065871543799 \bmod{11^{30}}\\
 \RL_{11}(\Hh,\Hf,\Hh)(5,2,5) & \equiv & 4791560577275108790581414445515 \bmod{11^{30}}.
\end{array}
\]
Using the Darmon-Rotger formula we compute
\[ \AJ_{11}(\Delta)(\eta_g^{{\rm u-r}} \otimes \omega_f \otimes \omega_g) \equiv -646073276230754578213318125190 \times 11 \bmod{11^{30}}.\]
Rather than attempt to recover a point from this, we take the generator $P = (0,0)$ for $E_f(\Q)$ and compute
$\log_{E_f}(P)$ and then try to determine a relationship. One finds
\[ \AJ_{11}(\Delta)(\eta_g^{{\rm u-r}} \otimes \omega_f \otimes \omega_g) \equiv \frac{258}{107} \log_{E_f}(P) \bmod{11^{30}}.\]
(We checked that multiplying the relevant projection operator by $2 \times 107$ gives an element in the integral
Hecke algebra.) Similarly we found
\[ \AJ_{11}(\Delta)(\eta_h^{{\rm u-r}} \otimes \omega_f \otimes \omega_h) \equiv -\frac{6708}{5647} \log_{E_f}(P) \bmod{11^{30}}.\]
\end{example}

The examples above suggest the construction in \cite{DRS} can be generalised, at least in a $p$-adic setting.

We now look at some examples in which one removes one of the main conditions in the Darmon-Rotger theorem \cite[Theorem 1.3]{DR-GrossZagier} itself, that the
prime does not divide the level. In each example rather than try to recover a rational point, we look for an algebraic
relationship between the logarithm of a generator and the special value we compute.

\begin{example}\label{Ex-77a}
Let $E_f\,:\, y^2 + y = x^3 + 2 x$
be the rank $1$ curve of conductor $77$ with Cremona label ``77a'' associated to the cusp form
\[ f :=  q - 3 q^3 - 2 q^4 - q^5 - q^7 + 6 q^9 - q^{11} + \cdots .\]
Let $g$ be the level $11$ and weight $2$ newform (associated to a rank zero elliptic curve):
\[ g := q - 2 q^2 - q^3 + 2 q^4 + q^5 + 2 q^6 - 2 q^7 - 2 q^9 - 2 q^{10} + q^{11} + \cdots .\]
We take the prime $p := 7$, which divides the level of $f$, and writing $\Hf$ and $\Hg$ for the Hida
families compute
\[ \RL_{7}(\Hg,\Hf,\Hg)(2,2,2)  \equiv   -1861584104004734313229493 \times 7 \bmod{7^{31}}.\]
Taking the generator $P = (2,3)$ we compute $\log_{E_f}(P) \bmod{7^{31}}$ and find that
$\log_{E_f}(P)/7\RL_{7}(\Hg,\Hf,\Hg)(2,2,2)$ satisfies the quadratic equation
$1600 t^2 + 48 t + 9 = 0$ modulo $7^{29}$.
\end{example}

The factor $p$ which occurs in the expression relating the special value to the logarithm of a point when the
prime divides the level is also seen in the next examples.

\begin{example}\label{Ex-469}
Let $E_f\,:\, y^2 + xy + y = x^3 - 80x - 275$ and $E_g: y^2 + xy + y = x^3 - x^2 - 12x + 18$
be the rank $1$ curves of conductor $469$ with Cremona labels ``469a'' and ``469b'', respectively, associated to the cusp forms
\[
\begin{array}{rcl}
f & := & q + q^2 + q^3 - q^4 - 3 q^5 + q^6 - q^7 - 3 q^8 - 2 q^9 - 3 q^{10} + \cdots\\
g & := & q - q^2 - 3 q^3 - q^4 + q^5 + 3 q^6 - q^7 + 3 q^8 + 6 q^9 - q^{10} + \cdots .
\end{array}
\]
Taking the prime $p := 7$ we compute
\[
\begin{array}{rcl}
\RL_{7}(\Hg,\Hf,\Hg)(2,2,2) & \equiv & 1435409545849510941783817 \bmod{7^{30}}\\
\RL_{7}(\Hf,\Hg,\Hf)(2,2,2) & \equiv &  6915472639041460159095363 \bmod{7^{30}}.
\end{array}
\]
Using generators $P_f = (-5,4)$ and $P_g = (2,-1)$ for $E_f(\Q)$ and $E_g(\Q)$, respectively, we found
\[
\begin{array}{rcl}
7 \RL_{7}(\Hg,\Hf,\Hg)(2,2,2)  & \equiv & 4 \log_{E_f}(P_f) \bmod{7^{30}}\\
35 \RL_{7}(\Hf,\Hg,\Hf)(2,2,2)  & \equiv & -16 \log_{E_g}(P_g) \bmod{7^{30}}.
\end{array}
\]
\end{example}

In the above example the ``tame'' level used in our computation was $N = 67 = \frac{469}{7}$. In the next example
it is one: for tame level one the author's algorithm does not use the theory of modular symbols at all, cf. \cite[Section 3.2]{AL}.

\begin{example}
Let $E_f\,:\, y^2 + xy + y = x^3 + x^2 - x $ 
be the rank $1$ curve of conductor $89$ with Cremona label ``89a'' associated to the cusp form
\[ f := q - q^2 - q^3 - q^4 - q^5 + q^6 - 4 q^7 + 3 q^8 - 2 q^9 + q^{10} - 2 q^{11} + \cdots .\]
Let $g$ be the level $89$ and weight $2$ newform (associated to a rank zero elliptic curve):
\[ g := q + q^2 + 2 q^3 - q^4 - 2 q^5 + 2 q^6 + 2 q^7 - 3 q^8 + q^9 - 2 q^{10} - 4 q^{11}  + \cdots .\]
Taking the prime $p := 89$ we found that
\[ 89 \RL_{89}(\Hg,\Hf,\Hg)(2,2,2) \equiv 72 \log_{E_f}(P) \bmod{89^{21}}\]
where $P = (0,0)$ is a generator.
\end{example}

The author understands that the above examples are consistent with on-going work of Darmon and Rotger to generalise
their formula to the situation in which the prime $p$ does divide the level $N$ \cite{DR-2}.

\end{document}